\documentclass[centertags]{article}

\usepackage{amsmath}
\usepackage{amsthm}
\usepackage{amscd}
\usepackage{amssymb}
\usepackage{verbatim}

\title{A H\"older continuous vector field tangent to many foliations}
\author{Christian Bonatti\thanks{Partially supported by Northwestern
University and the Institut de Math\'ematiques de Bourgogne.  I thank
the Northwestern Math. Department for its warm hospitality during
visits while this paper was prepared.} \ and John
Franks\thanks{Supported in part by NSF grant DMS9803346.}}

\newtheorem{theorem}{Theorem}[section]

\newtheorem{quest}[theorem]{Question}
\newtheorem{exmple}[theorem]{Example}

\newtheorem{thm}{Theorem}[section] 

\newtheorem{lem}[thm]{Lemma} 
\newtheorem{prop}[thm]{Proposition}
 
\newtheorem{rem}[thm]{Remark} 

\newcommand{\R}{{\mathbb R}}
\newcommand{\Z}{\mathbb Z}

\newcommand{\F}{\mathcal F}
\newcommand{\C}{\mathcal C}
\newcommand{\N}{\mathbb N}

\begin{document}
\maketitle
\begin{abstract}
We construct an example of a H\"older continuous vector field on the
plane which is tangent to all foliations in a continuous family of
pairwise distinct $C^1$ foliations.  Given any $1 \le r <\infty,$ the
construction can be done in such a way that each leaf of each
foliation is the graph of a $C^r$ function from $\R$ to $\R.$
We also show the existence of a continuous vector field $X$ on $\R^2$
and two foliations $\cal{F}$ and $\cal{G}$ on $\R^2$ each
tangent to $X$ with a dense subset $\cal E$ of $\R^2$ such that at
every point $x\in \cal E$ the leaves $F_x$ and $G_x$ of the foliation
$\cal{F}$ and $\cal{G}$ through $x$ are topologically transverse.
\end{abstract}

\section{Introduction}

A basic result of ordinary differential equations asserts that any
non-vanishing Lipschitz vector field has a unique integral curve
through any point.  Simple examples show that this is not true when
the Lipschitz property is weakened to a H\"older condition.  For
example the differential equation $y' = \sqrt{|y|}$ with initial
condition $y(0) = 0$ admits the solutions $y(x) = x^2$ and $y(x)$
identically zero. A classical consequence of this non-unique
integrability is the existence of H\"older continuous vector fields on
the plane $\R^2$ which are not tangent to any foliation: for
example the vector field $\frac\partial{\partial x}
+\varphi(y)\frac\partial{\partial y}$ where $\varphi(y)=\sqrt y$\ if
$y\geq 0$ and $\varphi(y)=0$ if $y\leq 0$.

We present here a somewhat more sophisticated construction: continuous
vector fields (without singularities) on $\R^2$ which are tangent to
many different foliations.

\begin{theorem}
\label{thm:example}
For any $ 1\le r < \infty$ there is a H\"older continuous vector
field $X$ on $\R^2$ with the property that there is
family of pairwise distinct $C^1$ foliations $\F_t$ for $t \in [0,1]$ such that
$X$ is tangent to each foliation.  Moreover, each leaf of each foliation is
the graph of a  $C^r$ function from $\R$ to $\R.$
\end{theorem}
In the example built for Theorem~\ref{thm:example}, the vector field
$X$ is uniquely integrable on a dense open subset of $\R^2$.

\begin{theorem}\label{t.transverse}
There exists a continuous unit vector field $X$ on $\R^2$
and two foliations $\cal{F}$ and $\cal{G}$ on $\R^2$ each of them
tangent to $X$, and a dense subset $\cal E$ of $\R^2$ such that at
every point $x\in \cal E$ the leaves $F_x$ and $G_x$ of the foliation
$\cal{F}$ and $\cal{G}$ through $x$ are topologically transverse.
\end{theorem}

(In fact the vector field $X$ is also tangent to an uncountable family
of foliations as noted in Remark~\ref{r.many}).  Unfortunately we
know neither whether the vector field of Theorem~\ref{t.transverse}
can be H\"older continuous, nor if the leaves of the foliations can be
arbitrarily smooth.

This work is motivated by the invariant foliation appearing in the
study of smooth dynamical systems. For uniform hyperbolic systems, the
stable and unstable bundles are in general no more than H\"older
continuous. However, for dynamical reasons, they are always uniquely
integrable. Moreover, they are absolutely continuous with respect
to Lebesgue measure see \cite{BR,S}. This property has been
generalized in the very general context of Pesin theory see
\cite{PS,PS}.

More recently, an important effort has been made to understand a large
class of systems which present a weak form of hyperbolicity, namely,
partially hyperbolic systems, see \cite{BP,HPS}. These systems present
invariant strong stable and strong unstable directions, but also have
an invariant direction on which the dynamics may have weak expansion
or contraction. Each of these bundles are H\" older
continuous. Dynamical arguments prove the existence of foliations
tangent to the strong stable and the strong unstable directions, and
\cite{PS} ensures that these foliations are absolutely continuous with
respect to Lebesgue measure.  The third invariant direction, called the
{\em central direction} (or {\em central bundle}) is still not
understood.  We conjecture that, if the central direction has dimension
greater that $2$, it may be non integrable. If the central direction
has dimension $1$, we do not know either if it is uniquely
integrable, or if there exists a foliation tangent to it, or if
there is an invariant foliation tangent to it (if there are many
foliations tangent to it, maybe none of them are invariant). Recently,
Shub and Wilkinson exhibited the first known example where the central
direction is tangent to a foliation which is not absolutely continuous
with respect to Lebesgue measure (see \cite{SW}).

This lack of knowledge on the central direction is an important problem
in this theory.  In fact when the central foliation is regular in a
certain sense (see \cite{HPS} for the precise statement) it is {\em
structurally stable}, i.e., each small $C^1$ perturbation of the
dynamics admits a central foliation which is conjugate to the initial
one.  Most of the known examples of {\em robustly transitive} (see
\cite{Sh,M,BD,BV}) diffeomorphisms or of {\em stably ergodic} see
\cite{GPS,BPSW} diffeomorphisms are based on this property.  Many
works on these classes of diffeomorphisms assume the existence of an
invariant central foliation.

\begin{quest} 
Do there exist robustly transitive or stably ergodic diffeomorphisms
of a closed $3-$manifold, having a $1-$dimensional central bundle
which is not tangent to a unique invariant foliation?
\end{quest}

In order to pursue this problem, Wilkinson asked recently if
there exists a continuous non-singular vector field of $\R^2$ tangent to more
than one foliation.  This article provides a positive answer to this
question. As the central bundles are always H\"older continuous, we try
also to give a H\"older continuous example, as in
Theorem~\ref{thm:example}.

\section{The construction for Theorem~\ref{thm:example}}

We begin by constructing a single leaf as the graph of a function
$g: \R \to \R.$  We specify it by giving its derivative $h: \R \to \R.$

Consider the ``middle third'' Cantor set $\C$ in $[0,1]$ obtained by removing
at ``stage $n$'' $2^n$ gaps of length $1/3^n$, the so-called middle thirds.
For $x \in \C$ we define $h(x) = 0.$  And if $x \in (a,b)$ where $(a,b) = (k/3^n, (k+1)/3^n)$ is one
of the complementary gaps of length $1/3^n$ for $\C$ we define 
\[
h(x) = 3^{-2rn}\phi(3^n(x-a))
\]
where $\phi(t) = t^{r+1}(1-t)^{r+1}$, a smooth non-negative function on
$[0,1]$ which vanishes to order $r$ at both endpoints.  For later use
we note that if $A = \int_0^1 \phi(t) dt$ then by a simple change of variables
\begin{equation}
\int_a^b h(t) dt = 3^{-3rn} A \label{eqn1}.
\end{equation}

We also note that if $K$ is an upper bound for $\phi$ and its first
$r$ derivatives on $[0,1]$ then $h(x)$ and its first $r$ derivatives
are bounded on $[a,b]$ by $3^{rn} K 3^{-2rn} = K 3^{-rn}$ which tends uniformly
to zero as $n$ tends to infinity.  It follows that $h$ is a $C^r$
function on $[0,1]$ which vanishes along with its first $r$
derivatives at points of $\C.$

We now define 
\[
g(x) = \int_0^x h(t) dt \text{ for }x \in [0,1].
\]
Note that $g$ is strictly monotonic increasing since $g(x_2) - g(x_1)
= \int_{x_1}^{x_2} h(t)dt$ and $h$ is non-negative and strictly
positive on a dense set.  Hence $g:[0,1] \to[0, g(1)]$ is a $C^{r+1}$
homeomorphism.  We may extend it to a $C^{r+1}$ homeomorphism of $\R$
by, for example, setting $g(x) = -x^{2r}$ for $x <0$ and $g(x) = g(1) +
(x-1)^{2r}$ for $x > 1.$

The graph of $g$ will be one leaf of one of the foliations we
construct.  The other leaves of this foliation will be horizontal
translates of this leaf.  More precisely, let $L_0$ be the graph of
$g$ and let $\F_0$ be the foliation of $\R^2$ whose leaves are the
curves $L_c, c\in \R$ where
\[
L_c = L_0 +(c,0) = \{(x, g(x - c)) | x \in \R\} = \{(g^{-1}(y) + c, y)) | y \in \R\}.
\]
We will define the vector field $X$ to be tangent to this foliation.  We
note that if $(x_0,y_0) \in L_c$ then $c = x_0 - g^{-1}(y_0)$.  The slope
of $L_c$ at $(x_0,y_0)$ is $h(x_0 - c) =  h(g^{-1}(y_0)).$  Hence we
may define
\[
X = \frac{\partial}{\partial x} + h(g^{-1}(y))\frac{\partial}{\partial y},
\]
which is clearly a continuous vector field since $h$ is $C^r$ and
$g^{-1}(y)$ is continuous.  We now want to show that $h(g^{-1}(y))$ is
H\"older continuous.

Clearly it suffices to show that there is are constants $C>0$ and $\alpha \in (0,1)$
such that for every $0 \le y_1 < y_2 \le g(1)$ 
\[
\frac{|h(g^{-1}(y_2)) -  h(g^{-1}(y_1))|}{|y_2 -y_1|^\alpha} \le C.
\]
Letting $x_1 = g^{-1}(y_1)$ and $x_2 = g^{-1}(y_2)$ we want
\[
\frac{|h(x_2) -  h(x_1)|}{|g(x_2) -g(x_1)|^\alpha} \le C,
\]
for every $0 \le x_1 < x_2 \le 1$ or equivalently we want to show
$|g(x_2) -g(x_1)| \ge C^{-1} |h(x_2) -  h(x_1)|^\beta$
for some $C > 0$ and some $\beta = 1/\alpha > 1.$

Since $h$ is Lipschitz on $[0,1]$ there is a constant $D >0$ such that
$|h(x_2) -  h(x_1)| \le D |x_2 -  x_1|.$  Hence it will suffice to find
$\beta$ and $K$ such that 
\begin{equation}
|g(x_2) -g(x_1)| \ge K|x_2 -  x_1|^\beta \label{eqn2}
\end{equation}
since then $|g(x_2) -g(x_1)| \ge KD^{-\beta}|h(x_2) -  h(x_1)|^\beta .$

To show this we let $n$ be the unique positive integer such that
$2/3^{n-1} \le |x_2 -  x_1| < 2/3^{n-2}.$   It follows that there is
an interval $[j/3^{n-1}, (j+1)/3^{n-1}] \subset [x_1,x_2]$.  If this is not one
of the Cantor set gaps then it contains one $[k/3^{n}, (k+1)/3^{n}]$.
Hence there is a Cantor set gap $[a,b] \subset [x_1,x_2]$ of length
at least $1/3^{n}.$  Therefore
\[
g(x_2) -g(x_1) = \int_{x_1}^{x_2} h(t) dt \ge \int_{a}^{b} h(t) dt \ge A 3^{-3rn}
\]
by equation (\ref{eqn1}) above.  Thus if $\beta = 3r$ 
\[
|g(x_2) -g(x_1)|\ge A 3^{-3rn} = A (3^{-n})^\beta \ge A (9/2)^\beta |x_2 -  x_1|^\beta,
\]
since $|x_2 - x_1| < 2/3^{n-2}.$ Hence we have established equation
(\ref{eqn2}) and the vector field $X$ is H\"older continuous.

We are now prepared to construct the other foliations tangent to the
vector field $X$.  Let $\C' = g(\C)$ so $\C'$ is a Cantor set in
$[0,g(1)].$ Let $\psi: [0,g(1)] \to \R$ be a Cantor function
associated with $\C'$.  More precisely $\psi(y)$ is a monotonic
increasing continuous function which is not constant, but is constant
on each component of $[0,g(1)] \setminus \C'$.  The construction of
such a function can be found in Royden, \cite{Ro}, page 39.  We extend
$\psi$ to all of $\R$ by setting $\psi(y) = \psi(0)$ for $y < 0$ and
$\psi(y) = \psi(g(1))$ for $y > 0.$

For each $t \in [0,1]$ we can consider the function $f_t: \R \to \R$
given by $x = f_t(y) = g^{-1}(y) + t\psi(y).$ This is a strictly
increasing function and easily seen to be a homeomorphism of $\R.$
Hence we may consider the inverse homeomorphism $g_t = f_t^{-1}.$
Denote by $\C_t$ the image $\C_t=f_t(\C')=f_t(g(\C))$. It is a Cantor
subset of $\R$.
 
We
wish to show that $g_t(x)$ is a $C^r$ function of $x$ and that its
graph is tangent to the vector field $X$.  This is easily seen to be
true on a neighborhood of any point $x_0$ which is not in the Cantor
set $\C_t$. This
is because the function $\psi(y)$ is constant on a
neighborhood of $y_0 = g_t(x_0)$ and hence if $t\psi(y) = c_0$ on such
a neighborhood then the curve $x = f_t(y) = g^{-1}(y) + t\psi(y)$ is
$x = g^{-1}(y) + c_0$ or $g(x - c_0) = y.$ This clearly implies that
$g_t$ is $C^r$ and tangent to $X$ on a neighborhood of $x_0.$

It remains to show that for fixed $t$ the function $g_t(x)$ is $C^r$
and tangent to $X$ at points of $\C_t.$ We do this by showing that $g_t$
is ``flatter'' than $g$ near points of $\C_t.$ 

Let $z_0 \in \C_t$ so $z_0 = f_t(g(x_0))$ for some $x_0 \in \C.$  If
$z \in \R$ then $z = f_t(g(x))$ for some $x.$  Note that
\[
z_0 = f_t(g(x_0)) = g^{-1}(g(x_0)) + t \psi(g(x_0)) = x_0 + t \psi(g(x_0))
\]
and similarly
\[
z = f_t(g(x)) = g^{-1}(g(x)) + t \psi(g(x)) = x + t \psi(g(x))
\]
so $z - z_0 = (x - x_0) + t (\psi(g(x)) - \psi(g(x_0)).$  Since
both $g$ and $\psi$ are monotonic increasing, we conclude that
\begin{equation}
|x-x_0| \le |z - z_0|. \label{eqn3}
\end{equation}

Then we observe that 
\begin{equation}
|g_t(z) - g_t(z_0)| = |g_t(f_t(g(x))) -
g_t(f_t(g(x_0)))| = |g(x) - g(x_0)|. \label{eqn4}
\end{equation}
But since $g$ is $C^{r+1}$ and
its first $r$ derivatives vanish at points of $\C$ we know there is a
constant $B$ such that $|g(x) - g(x_0)| \le B |x -x_0|^{r+1}.$
Combining this with equations (\ref{eqn3}) and (\ref{eqn4}) 
above we conclude that
\[
|g_t(z) - g_t(z_0)| \le B|z - z_0|^{r+1}.
\]
This implies that $g_t$ is ${r}$ times differentiable at $z_0$ and
that its first $r$ derivatives vanish there.  In particular the graph
of $g_t$ is tangent to the vector field $X$.

We define the foliation $\F_t$ to have leaves which are horizontal
translates of the graph of $g_t$.  That is, we let $L^t_0$ be the
graph of $g_t$ and define $L^t_c = L^t_0 +(c,0)$.  The leaf $L^t_c$ is
the graph of $y = g_t(x -c)$.  Since the vector field $X$ is invariant
under horizontal translation, we conclude that $L^t_c$ is tangent to
$X.$ Hence for each fixed $t \in [0,1]$ the foliation $\F_t$ has $C^r$
leaves all of which are tangent to the vector field $X$.

\section{The construction for Theorem~\ref{t.transverse}}

Our construction is based on the following remark:

\begin{rem} 
Let $X$ be a $C^0$ non-singular vector field on $\R^2$. Assume
that there is a family $\{\gamma_i\}_{i\in \N}$ of proper embeddings
$\gamma_i : \R\to \R^2$ with the following 3
properties:
\begin{enumerate}
\item the curves $\gamma_i$ are pairwise disjoint,
\item  the union of the curves $\gamma_i$ is dense in $\R^2$,
\item each of the $\gamma_i$ is everywhere tangent to $X$.
\end{enumerate}

Then there exists a unique foliation $\cal{F}$ admitting the $\gamma_i$ as leaves. Moreover $\cal{F}$ is everywhere tangent to $X$. 
\label{r.foliation}
\end{rem}

We will obtain the announced example as a limit of a sequence of
foliations ${\cal{F}}_n$ and construct a family $\{\gamma_n\}$ of
curves such that each curve $\gamma_n$ is tangent to, but
topologically transverse to, the foliations ${\cal{F}}_m$, $m>n$.
Moreover, the curves $\gamma_n$ are pairwise disjoint and their union
is dense in $\R^2$.  We will show that the foliations
${\cal{F}}_n$ converge to a foliation $\cal{F}$ whose leaves are
tangent to, but topologically transverse to, each curve
$\gamma_n$. Then the family of curves $\{\gamma_n\}$ will be completed to
a foliation $\cal{G}$ and we will show that the leaves of $\cal{G}$
and $\cal{F}$ through each point $x$ are tangent to the same
continuous vector field. This construction is summarized by the
following proposition:

\begin{prop} 
Let $E_1,E_2$ denote the unit vector fields parallel to the canonical
basis of $\R^2$. There is a sequence
$(Z_n,{\cal{F}}_n,\gamma_n)_{n\in \N}$ having the following
properties:
\begin{enumerate}
\item \label{i.Zn} The sequence $\{Z_n\}$ of continuous
unit vector fields converges uniformly to a vector field
$Z$. Moreover we may assume that the coordinates of $Z_n$ in the $E_1,E_2$
basis belong to $[1/3,2/3]$.
\item \label{i.Fn} The foliation $\F_n$ is tangent to $Z_n$. 
\item \label{i.gamman} The $\gamma_n$ are proper embeddings of $\R$
into $\R^2$ pairwise disjoint and tangent to any of the vector fields
$Z_m$, $m>n$, but transverse (positively) to $Z_n$, and the union
$\bigcup_n\gamma_n$ is dense in $\R^2$,
\item \label{i.transverse} The $\gamma_n$ cut topologically
transversally and positively each leaf of each $\F_m$, $m\geq n$, 
in exactly one point: thus the holonomy map $\varphi_{m,i,j}$ of the
foliation ${\cal{F}}_m$ from $\gamma_i$ to $\gamma_j$ is well defined
for $i\leq m$ and $j\leq m$
\item \label{i.holonomy} The holonomy map $\varphi_{m,i,j}$ does not
depend on $m$, that is, $\varphi_{m,i,j}=\varphi_{m+1,i,j}$.
\item For any $n>0$ and any $x\in \gamma_0$ we denote by $F_{n,x}$ the leaf
of ${\cal F}_n$ through $x$. There is a sequence
$\varepsilon_n>0$ with the following properties:
\begin{enumerate}
\item \label{i.a} for any point $y\in F_{n,x}$, the horizontal
projection $y'$ of $y$ on the leaf $F_{n+1,x}$ satisfies
$d(y,y')<\varepsilon_n$
\item \label{i.b}
for any $x_1,x_2\in \gamma_0$ with $\|x_1\|\leq n$ and $d(x_1,x_2)\geq \frac1n$ , the distance $d(F_{n,x_1},F_{n,x_2})$ is greater than $10\cdot\sum_{i\geq n}\varepsilon_i$; that is, 
$$
\inf\{d(y_1,y_2), y_1\in F_{n,x_1} \mbox{ and } y_2\in F_{n,x_2}\}\geq 10\cdot\sum_{i\geq n}\varepsilon_i
$$
\end{enumerate}
\end{enumerate}
\label{p.transverse}
\end{prop}

Before proving this proposition, we show that it proves
Theorem~\ref{t.transverse}.
\vskip 2mm

\noindent
{\bf Proof of Theorem~\ref{t.transverse}:} First notice that
item~\ref{i.gamman} implies that the curves $\gamma_n$ are all tangent
to $Z$, and satisfy all the hypotheses of Remark~\ref{r.foliation}, so
that the family $\gamma_n$ can be completed in a unique way to a
foliation $\cal G$ tangent to $Z$.

{From} item 1 of Proposition~\ref{p.transverse} the leaf $F_{n,x}$ of
${\cal F}_n$ through $x\in\gamma_0$ can be seen as the graph of a function
$f_{n,x}\colon\{0\}\times \R\to \R\times\{0\}$, and
item~\ref{i.a} and \ref{i.b} imply that for any $x\in\gamma_0$ the
functions $f_{n,x}$ converge uniformly to some function $f_x$.
Item 1 implies that the graph $F_x$ of $f_x$ is tangent to the limit
vector field $Z$.

Item~\ref{i.a} and \ref{i.b} imply that for $x_1\neq x_2$ the curves
$F_{x_1}$ and $F_{x_2}$ are disjoint as we now show.
Choose an integer $n$ such
that $\|x_1\|\leq n$ and $d(x_1,x_2)\geq \frac1n$; by item~\ref{i.a},
for any point $y_1\in F_{n,x_1}$ and $y_2\in F_{n,x_2}$, the
horizontal projections $y'_i$ of $y_i$ on the leaf $F_{x_i}$, with
$i=1, 2$, satisfy $d(y_i,y')<\sum_{j\geq n}\varepsilon_j$. Then
\[
\inf\{d(y'_1,y'_2), y'_1\in F_{x_1} \text{ and } y'_2\in F_{x_2}\}
\]
is greater than
\[
\inf\{d(y_1,y_2), y_1\in F_{n,x_1} \text{ and }
y_2\in F_{n,x_2}\}- 2\cdot\sum_{j\geq n}\varepsilon_j, 
\]
and
item~\ref{i.b} implies that this distance is greater than
$8\cdot\sum_{j\geq n}\varepsilon_j$.

Each $\F_n$ is a foliation, so the union of the
$F_{n,x}$ are all of $\R^2$, for every $n$. Since, for $n$ large the leaves of
$F_x$ are (uniformly in $x$ and $n$) close to the leaves $F_{n,x}$, we
get that the union of the $F_x$ is dense in $\R^2$, so that the
$F_x$ are leaves of a unique foliation $\cal F$, tangent to $Z$.

To finish the proof of Theorem~\ref{t.transverse} it remains to note
that the foliations $\cal F$ and $\cal G$ are topologically transverse
along each of the $\gamma_n$, that is, each curve $\gamma_n$ cuts
each leaf of $\cal F$ in exactly one point. This is a direct
consequence of the fact that the holonomy map $\varphi_{i,j}$ of $\cal
F$ from $\gamma_i$ to $\gamma_j$ is well defined and coincides with
$\varphi_{m,i,j}$, for $m>\sup\{i,j\}$.

\hfill$\square$
\medskip
\begin{rem}\label{r.many}
In fact the vector field $Z$ constructed in Theorem~\ref{t.transverse}
above is tangent to many foliations: consider the vertical strips
bounded by the vertical lines $\{n\}\times \R$ where $n\in\Z$ (which
are transverse to $Z$, by construction). Then, to any point
$\omega\in\{F,G\}^{\Z}$ we associate the foliation ${\cal F}_\omega$
whose leaves coincide in $[n,n+1]\times \R$ with the foliation $\cal
F$ or $\cal G$ according with the $n^{th}$ letter of the infinite word
$\omega$. One verifies easily that ${\cal F}_\omega$ is a foliation
tangent to $Z$.
\end{rem}

The proof of Proposition~\ref{p.transverse} follows directly from the
following lemmas:

\begin{lem} 
Let $X$ be a continuous vector field on $\R^2$, and let $\gamma$ be a
proper embedding of $\R$ in $\R^2$, transverse to $X$.  Assume that
there is a neighborhood $U$ of $\gamma$ on which $X$ is uniquely
integrable.

Then there is a tubular neighborhood $V\subset U$ of $\gamma$ which is
the image of a proper embedding $\psi\colon [0,1]\times\R\to \R^1$
such that $\gamma_0=\psi(\{0\}\times \R)$ and
$\gamma_1=\psi(\{0\}\times \R)$ are transverse to $X$ and each orbit
of $X$ in $V$ is a segment joining $\gamma_0$ to $\gamma_1$, so that
the holonomy is a homeomorphism $\varphi\colon\gamma_0\to \gamma_1$.

\label{l.tubular}
\end{lem}
The proof of this lemma is a simple exercise, using the local flow of
$X$ defined on $U$, where it is smooth and so uniquely integrable.

\begin{lem} 
Under the assumptions of Lemma~\ref{l.tubular} and with the same
notation, assume that there is a $\delta>0$ and a continuous function
$h:\gamma\to ]0,\delta]$ such that the tangent direction to
$\gamma$ at any point $x\in \gamma$ is given by $X(x)+h(x)Y(x)$
where $Y(x)$ is a unit vector orthogonal to $X(x)$.

Then there is a continuous vector field $\tilde X$ satisfying
the following properties: 
\begin{enumerate}
\item $\tilde X$ coincides with $X$ outside of the tubular neighborhood $V$
\item $\tilde X$ is smooth on $V\setminus\gamma$,
\item $\tilde X$ is tangent to $\gamma$ at every point $x\in\gamma$,
\item The vector field $\tilde X$ is tangent in $V$ to a unique
foliation $\tilde{\cal F}$ whose leaves are segments joining $\gamma_0$
to $\gamma_1$, and the holonomy from $\gamma_0\to \gamma_1$ is
precisely $\varphi$ (that is coincides with the holonomy associated to
$X$),
\item The vector field $\tilde X$ is a $2\delta$-$C^0$-perturbation of $X$ that is,  
$\|X(x)-\tilde X(x)\|\leq 2\delta$, for any $x\in \R^2$.
\item 
For any point $x$ outside of the tubular neighborhood $V$, the leaf
$\tilde F_x$ of $\tilde{\F}$ through $x$ is $2\delta-$close to
the leaf $F_x$ of $\cal F$ through $x$.
\end{enumerate}
\label{l.local}
\end{lem}
\noindent
{\bf Proof:} We choose a very small closed tubular neighborhood
$W\subset V$ of $\gamma$ bounded by two curves $\lambda_0$ and
$\lambda_1$ transverse to $X$, such that $\gamma_0\cup\lambda_0$ and
$\lambda_1\cup\gamma_1$ are the boundary of two strips $W_0$ and
$W_1$, respectively, trivially foliated by $\cal F$, which are the
closure of the connected components of $V\setminus W$.

We first perturb $X$ in the interior of $W$ in order to get a
topological vector field $\tilde X$ in $W$ tangent to $\gamma$ and
tangent to a unique foliation $\tilde{\cal F}$ (topologically
transverse to $\gamma$) whose leaves are segments joining $\lambda_0$
to $\lambda_1$, without taking care of the holonomy from $\lambda_0$
to $\lambda_1$; if $W$ has been chosen small enough, each segment of
leaf of $\tilde {\cal F}$ joining a point of $\lambda_0$ to
$\lambda_1$ has a length very small in comparison with $\delta$.

Now we perturb $X$ in the interior of the strips $W_0$ and $W_1$ in
order to recover the property that the holonomy of $\tilde{\cal F}$
from $\gamma_0$ to $\gamma_1$ (i.e. the composition of the holonomies from
$\gamma_0$ to $\lambda_0$, from $\lambda_0$ to $\lambda_1$ and finally
from $\lambda_1$ to $\gamma_1$) coincides with the holonomy of $\cal
F$.  \hfill$\square$
\medskip

The main technical difficulty for proving Proposition~\ref{p.transverse} is:
\begin{lem} Let $X$ be a continuous non-singular vector field of $\R^2$, 
tangent to a foliation $\cal F$ conjugate to the trivial
foliation. Assume that there is a family
$\sigma_1,\sigma_2,\gamma_1,\dots,\gamma_k$ of proper smooth
embeddings of $\R$ to $\R^2$ with the following properties:

\begin{enumerate}
\item $\sigma_1$ and $\sigma_2$ are transverse to $X$ and cut each
leaf of $\cal F$ in a point.
\item $\sigma_1$ and $\sigma_2$ are disjoint so that
$\R^2\setminus (\sigma_1\cup \sigma_2)$ has three connected components,
$P_1$, $P_2$ and $P_3$ with $\partial(P_1)=\sigma_1$, $\partial
(P_3)=\sigma_2$ and $\partial(P_2)=\sigma_1\cup \sigma_2$. Each of the 
$\gamma_i$ lies in $P_2.$ One may assume
that $X$ coincides on $P_1\cup P_3$ with the constant vector field
$\frac\partial{\partial x}+\frac\partial{\partial y}$. 
\item the curves $\gamma_i$ are pairwise disjoint, and are tangent to $X$.
\item the curves $\gamma_i$ are topologically transverse to $\cal F$,
and they cut the leaves of $\cal F$ with positive orientation
($\gamma_i$ and the leaf of $\cal F$ are oriented by $X$).
\item The vector field $X$ is smooth on $\R^2\setminus\bigcup_i\gamma_i$.
\end{enumerate} 

Then given any point $x_0\in \R^2\setminus\bigcup_i\gamma_i$ and any
$\eta>0$, there is a proper embedding
$\gamma_{k+1}\colon\R\to\R^2$ with $\gamma_{k+1}(0)=x_0$, disjoint from $\bigcup_i\gamma_i$,
transverse to $X$ and cutting each leaf of $\cal F$ positively in
exactly one point, and such that the tangent line at each point
$x\in\gamma_{k+1}$ is given by $X(x)+\psi(x)Y(x)$ where $Y(x)$ is a
unit vector orthogonal to $X$ and $0 <\psi(x)<\eta$.
\label{l.induction}
\end{lem}
\noindent
{\bf Proof:} Let $Y$ denote the unit vector field such that the basis
$(Y,X)$ is positively oriented and orthonormal.  One can choose
pairwise disjoint tubular neighborhoods $\Gamma_i$ of the curves
$\gamma_i$ such that their boundaries consist of smooth curves
$\gamma_i^+$ and $\gamma_i^-$ transverse to $X$, and such that
$\bigcup_i\Gamma_i$ does not contain the point $x_0$.  Moreover,
we can choose them so the tangent
to the curves $\gamma_i^\pm$ can be written $X+h_i^\pm(x)Y$ with
$0 < h_i^\pm(x)$. By
convention, the leaves of $\cal F$, oriented by $X$, enter $\Gamma_i$
through $\gamma_i^-$ and exit $\Gamma_i$ through $\gamma_i^+$.

Now we choose a smooth function $\psi(x)$ on $\R^2$, coinciding
with $\eta/2$ on a neighborhood of $\bigcup_i\gamma_i$ and on
$P_1\cup P_3$, such that $0<\psi(x)<\eta$ for any $x\in \R^2$
and such that
\begin{equation} \label{e.psi} 
0 < \psi(x) < h_i^\pm(x)
\end{equation} 
for any $x\in \gamma_i^\pm$. Let $Z=X+\psi(x)Y$, so it is a continuous
vector field which is smooth on $\R^2\setminus \bigcup_i\gamma_i$. We
claim that a maximal solution of this vector field through the point $x_0$
satisfies all the announced properties. First notice that as $Z$ is a
continuous vector field without singularity, any maximal solution is a
proper embedding of $\R$ in $\R^2$.

Consider a maximal solution $\gamma_+$ for the vector field $Z$ for
positive time and starting at $x_0$. We first show that it is disjoint
from all the curves $\gamma_i$.  The inequality~\ref{e.psi} implies
that $Z$ must enter $\Gamma_i$ through $\gamma_i^-$ and exit
$\Gamma_i$ through $\gamma_i^+$. So either it is disjoint from the
$\Gamma_i$ or it crosses some $\gamma_i^-$. Notice that it cannot cross
$\gamma_i^-$ twice. Moreover, it cannot contain a point
of $\gamma_i$.  This is because there is an infinite strip bounded by
$\gamma_i$ and $\gamma_i^-$ and the vector field $Z$ 
points inward on both boundary components of this strip.  Hence,
since $\gamma^+$ can only enter $\Gamma_i$ by crossing $\gamma_i^-$,
it cannot intersect $\gamma_i.$

The same argument shows that a maximal solution for negative time is
also disjoint from each $\gamma_i$.  Hence this solution remains in
$\R^2\setminus\bigcup_i \gamma_i$ where $Z$ is smooth, so that this
solution is unique and we can speak of the orbit of $x_0$ for $Z$

To finish the proof it remains to show that this orbit cuts all the
leaves of $\cal F$. {\it A priori} it cuts an interval of the space of
leaves of $\cal F$ (which is homeomorphic to $\R$ by hypothesis). Look
at the positive orbit $\gamma_+$ of $x_0$. If there is a sequence of
times $t_n\to+\infty$ such that $\gamma_+(t_n)\in P_2$ then the
corresponding interval has to be infinite, because the leaves of $\cal
F$ in $P_2$ are compact segments joining $\sigma_1$ to $\sigma_2$. On
the other hand, if for all large $t$ the point $\gamma_+(t)$ belongs to  $P_1 \cup P_3$,
then it follows from the fact that $\psi=\eta/2$ on
$P_1\cup P_3$ that the interval in the leaf space is infinite.  
A similar argument for the negative orbit through $x_0$ completes
the proof.
\hfill $\square$

\medskip

We are now ready to prove Proposition~\ref{p.transverse} finishing the proof of Theorem~\ref{t.transverse}:
\medskip

\noindent
{\bf Proof of Proposition~\ref{p.transverse}:} The sequence
$(Z_n,{\cal F}_n,\gamma_n,\varepsilon_n)$ is constructed beginning with the
trivial foliation ${\cal F}_0$ tangent to the constant vector field
$\frac\partial{\partial x}+\frac\partial{\partial y}$ and
inductively  using Lemmas~\ref{l.tubular}, ~\ref{l.local} and
\ref{l.induction}. The only point which is not straightforward is
the choice of the sequence $\varepsilon_n$, in order to get the
property of item~\ref{i.b}.

Fix a sequence of points $\{x_n\}$ which is dense in $\R^2$.

Assume that the foliation ${\cal F}_n$ has been built. By
construction all the leaves through a point $x\in \gamma_0$ with
$\|x\|<n$ coincide with a leaf of ${\cal F}_0$ outside of some compact
set. Thus (by continuity and compactness) there is some $\mu_n>0$
such that if $x,y\in \gamma_0$ satisfy $\|x\| \leq n$ and $d(x,y)\geq \frac1n$
the leaves $F_{n,x}$ and $F_{n,y}$ remain at a distance greater than
$\mu_n$. Then we choose
$\varepsilon_{n}<\frac1{10}\inf\{\varepsilon_{n-1},\frac12\mu_n\}$. This choice 
implies:
 $$\sum_{i\geq n}\varepsilon_i\leq \varepsilon_n\cdot\sum_0^\infty\frac1{10^i}<2\varepsilon_n<\frac1{10}\mu_n.$$

Choosing the constant $\delta$ in Lemma~\ref{l.local} sufficiently
small, any foliation obtained by perturbing ${\cal F}_n$ using this
lemma will satisfy item~\ref{i.a} with this choice of $\varepsilon_n$.  
If $x_{n+1}$ belongs to the union
of the $\gamma_i$, $i\leq n$ then we replace $x_{n+1}$ by the first
point of the sequence $\{x_i\}_{i > n}$ which is not in the union of
these curves.  Then Lemma~\ref{l.induction} allows us to choose a
curve $\gamma_{n+1}$ through $x_{n+1}$ disjoint from
$\bigcup_1^n\gamma_i$; we choose the constant $\eta$ in
Lemma~\ref{l.induction} less than $\delta$. Then we choose a
tubular neighborhood of $\gamma_{n+1}$ disjoint from
$\bigcup_1^n\gamma_i$ (using Lemma~\ref{l.tubular}) and we build
${\cal F}_{n+1}$ as in Lemma \ref{l.local}.

In order to apply  Lemma~\ref{l.induction} recursively, it remains
to get the curves $\sigma_{1,n+1}$ and $\sigma_{2,n+1}$. The existence
of these curves follows immediately from the facts that ${\cal F}_{n+1}$
coincides with ${\cal F}_n$ outside of the tubular neighborhood $V_n$ and
that the boundary of this neighborhood consists of two curves
transverse to ${\cal F}_n$ and cutting all the leaves of ${\cal F}_n$.

\vskip 3cm
\noindent
Christian Bonatti:\\
Institut de Math\'ematiques de Bourgogne, UMR 5584,\\ Universit\'e de Bourgogne\\
21 078 Dijon cedex, FRANCE\\
E-mail: bonatti@u-bourgogne.fr
\medskip

\noindent
John Franks:\\
Dept. of Mathematics, Northwestern University\\
Evanston, IL 60208, USA\\
E-mail: john@math.northwestern.edu

\end{document}